\documentclass{article}
\usepackage[american]{babel}
\usepackage{amsfonts,amsmath,amssymb,epsf,epsfig}
\usepackage{xypic}
\xyoption{all}
\input{xypic}
\newdir{ >}{{}*!/-8pt/\dir{>}}

\def\R{\mathbb{R}}

\def\Z{\mathbb{Z}}        
      
\def\T{\mathbb{T}}

\def\K{\mathbb{K}} 
    
\def\End{{\rm End}}
\def\coker{{\rm coker}}\def\ker{{\rm ker}}
\def\id{\rm id}\def\im{{\rm im}}

\def\pr{\noindent $\bf{Proof.}$\quad}     
\def\fin{\hfill$\square$\\}           
\newtheorem{theo}{Theorem}
\newtheorem{defi}{Definition}
\newtheorem{Rem}{{\bf Remark}}
\newenvironment{rem}{\begin{Rem} \strut\\ \normalfont}{\end{Rem}}
\newtheorem{prop}{Proposition}
\newtheorem{cor}{Corollary}

\newtheorem{lem}{Lemma}

\def\cP{{\mathcal P}}

\begin{document}
\title{On $2$-Holonomy} 

\author{Hossein Abbaspour\\
        Universit\'e de Nantes\\
\and    Friedrich Wagemann \\
        Universit\'e de Nantes}
   
\maketitle


\begin{abstract}
We construct a cycle in higher Hochschild homology
associated to the $2$-dimensional torus
which represents $2$-holonomy of a non-abelian gerbe
in the same way the ordinary holonomy of a principal $G$-bundle gives
rise to a cycle in ordinary Hochschild homology.
This is done using the connection $1$-form of Baez-Schreiber.
 
A crucial ingredient in our work 
is the possibility to arrange that in the structure crossed module
$\mu:{\mathfrak h}\to{\mathfrak g}$ of the principal $2$-bundle, the Lie algebra
${\mathfrak h}$ is abelian, up to equivalence of crossed modules.
\end{abstract}

\vspace{.5cm}



\tableofcontents

\section*{Introduction}

Principal $2$-bundles have been studied 
in \cite{BaeSch}, \cite{Bar}, \cite{GinSti}, \cite{MarPic}, \cite{SchWal1}, 
\cite{SchWal2}, \cite{Woc} and \cite{CLS} (This list is by no means 
exhaustive). 
We will sketch the different approaches and explain our point 
of view, namely, we choose a framework at the intersection of gerbe theory and 
theory of principal $2$-bundles. The structure $2$-group of a 
principal $2$-bundle is in our framework a strict $2$-group (we refrain
from considering more general structure groups like coherent $2$-groups), 
and its Lie algebra a strict Lie $2$-algebra, opening the way to using 
all information about strict Lie $2$-algebras which we discuss in the 
first section. 

The first (non-gerbal) approach to principal $2$-bundles is due to Bartels 
\cite{Bar}. He 
defines $2$-bundles by systematically categorifying spaces, groups and 
bundles. Bartels writes down the necessary coherence relations for a locally 
trivial principal $2$-bundle with structure group a coherent $2$-group. 
This work has then been taken up by Baez and 
Schreiber \cite{BaeSch} in order to define connections for principal 
$2$-bundles. In parallel work, Schreiber and 
Waldorf \cite{SchWal1}, \cite{SchWal2}, and Wockel \cite{Woc} also 
take up Bartels work in order to define holonomy (Schreiber-Waldorf)
or to pass to gauge groups (Wockel). Baez and 
Schreiber describe an approach using locally trivial $2$-fibrations whose
typical fiber is a strict $2$-group.   

Non-abelian gerbes and principal $2$-bundles are two notions which
are close, but have subtle differences. The cocycle data of the two notions
has been compared in \cite{BaeSch}, section 2.1.4 and 2.2. Baez and Schreiber 
show that under certain conditions, the description in terms of local data of
a principal $2$-bundle with $2$-connection is equivalent to the cocycle 
description of a (possibly twisted) non-abelian gerbe with {\it vanishing fake 
curvature}. 
This constraint is also shown to be sufficient for the existence
of $2$-holonomies, i.e. the parallel transport over surfaces.

The approach of Schreiber and Waldorf \cite{SchWal1}, \cite{SchWal2} is based 
on so-called transport functors. Schreiber and Waldorf push the equivalence 
between categories
of principal $G$-bundles with connection over $M$ and transport functors from 
the thin fundamental groupoid of $M$ to the classifying stack of $G$ to
categorical dimension $2$. These transport functors can then be described in 
terms of differential forms, i.e. for a trivial principal $G$-bundle, these 
transport functors correspond to $\Omega^1(M,{\mathfrak g})$, where 
${\mathfrak g}$ is the Lie algebra of $G$. They show similarly that 
$2$-transport functors from the thin fundamental $2$-groupoid correspond to
pairs of differential forms $A\in\Omega^1(M,{\mathfrak g})$ and 
$B\in\Omega^1(M,{\mathfrak h})$ with vanishing fake curvature $F_A+\mu(B)=0$, 
where $\mu:{\mathfrak h}\to
{\mathfrak g}$ is the crossed module of Lie algebras corresponding to the 
strict Lie $2$-group which comes into the problem. It is clear that this
approach is based on the notion of holonomy. 

Wockel \cite{Woc} also takes up Bartel's work. In order to make them more
easily accessible, he formulates a principal $2$-bundle over $M$ in terms of
spaces with a group action. A (semi-strict) principal $2$-bundle over $M$ is 
then a locally trivial ${\mathcal G}$-$2$-space. The $2$-group ${\mathcal G}$ 
is strict, and so is the action functor, but the local triviality requirement 
is not necessarily strict. Wockel shows that semi-strict principal $2$-bundles
over $M$ are classified by non-abelian \v{C}ech cohomology.

The approach of Ginot and Sti\'enon \cite{GinSti} is based on looking at a 
principal $G$-bundle as a generalized morphism (in the sense of Hilsum and 
Skandalis) from $M$ to $G$, both being considered as groupoids. In the same 
way they view principal $2$-bundles as generalized morphisms from the 
manifold $M$ (or in general some stack, represented by a Lie groupoid) 
to the $2$-group ${\mathcal G}$, both being viewed
as $2$-groupoids. In this context, they exhibit a link to gerbes (in their
incarnation as extensions of groupoids) and define characteristic classes.   

The particularity of Martins and Picken's approach \cite{MarPic}
is that they consider
special ${\mathcal G}$-$2$-bundles. For a strict $2$-group ${\mathcal G}$
whose associated crossed module is $\mu:H\to G$, these bundles are obtained 
from a principal $G$-bundle $P$ on $M$. The speciality requirement is  
that the principal ${\mathcal G}$-$2$-bundle is given by a non-abelian cocycle
$(g_{ij},h_{ijk})$ as below, but with $\mu(h_{ijk})=1$ in order 
to have a principal $G$-bundle $P$. Using the language which we will introduce 
below, Martins and Picken suppose that the band of the gerbe (which is in 
general a principal $G/\mu(H)$-bundle) lifts to a principal $G$-bundle.  
Martins and Picken define connections for these special 
${\mathcal G}$-$2$-bundles and $2$-holonomy $2$-functors. 

Chatterjee, Lahiri and Sengupta \cite{CLS} use in the first place a reference
connection $1$-form $\bar{A}$ in order to take for a fixed $G$-principal 
bundle $P\to M$ only $\bar{A}$-horizontal paths in the path space 
${\mathcal P}_{\bar{A}}P$ they consider. ${\mathcal P}_{\bar{A}}P$ is a 
$G$-principal bundle over the usual path space ${\mathcal P}M$.
Then given a pair $(A,B)$ as above,
they construct a connection $1$-form $\omega_{(A,B)}$ on ${\mathcal P}_{\bar{A}}P$
using Chen integrals. Major issues are reparametrization invariance and 
the curvature. The authors switch to a categorical description motivated by 
their differential geometric study in the end of the article.

Let us summarize the different approaches in the following table:

\hspace{.5cm}
\begin{tabular}{|c|l|} \hline
author(s) & concept \\ \hline\hline
Bartels & principal $2$-bundles with\\ 
& coherent structure group \\ \hline
Baez-Schreiber & global connection $1$-form \\
& for principal $2$-bundles\\ \hline
Schreiber-Waldorf & holonomy in terms of\\
& transport functors\\ \hline
Wockel & relation to non-abelian \\
& \v{C}ech cohomology \\  \hline
Ginot-Sti\'enon & $2$-bundles as Hilsum-Skandalis'\\ 
& generalized morphisms\\  \hline
Martin-Picken & connections and holonomy for \\
& {\it special} principal $2$-bundles\\  \hline
Chatterjee-Lahiri-Sengupta & connections and holonomy \\
& using $\bar{A}$-horizontal paths \\ 
& for a reference $1$-form $\bar{A}$ \\ \hline
\end{tabular}  
\vspace{.5cm}

The goal of our article is to construct a cycle in higher Hochschild homology
which represents $2$-holonomy of a non-abelian gerbe as described above
in the same way the ordinary holonomy gives
rise to a cycle in ordinary Hochschild homology, see \cite{AbbZei}.
This is done using the connection $1$-form of Baez-Schreiber \cite{BaeSch}
which we construct here from the band of the non-abelian gerbe.
 
A crucial ingredient in our work 
is the possibility to arrange that in an arbitrary crossed module
of Lie algebras $\mu:{\mathfrak h}\to{\mathfrak g}$, the Lie algebra
${\mathfrak h}$ is abelian, up to equivalence of crossed modules. This is 
shown in Section $1$ (see \cite{Wag}). The possibility to have ${\mathfrak h}$ 
abelian is used in order to obtain a {\it commutative} differential graded 
algebra $\Omega^*:=\Omega^*(M,U{\mathfrak h})$ whose higher Hochschild homology
$HH_{\bullet}^{\T}(\Omega^*,\Omega^*)$ associated to the $2$-dimensional torus 
$\T$ houses the holonomy cycle. We don't know of any definition of higher 
Hochschild homology for arbitrary differential graded algebras, therefore
we believe the reduction to abelian ${\mathfrak h}$ to be crucial when 
working with possibly non-abelian gerbes.
Section $1$ also provides a fundamental result on strict Lie $2$-algebras 
directly inspired from \cite{BaeCra}, namely, 
we explicitely show that the classification of strict Lie $2$-algebras
in terms of skeletal models (of the associated semi-strict Lie $2$-algebra)
and in terms of the associated crossed module coincide.     

Section $2$ reports on crossed modules of Lie groups. These play a minor role
in our study, because the main ingredient for the connection data is the 
the infinitesimal crossed module, i.e. the Lie algebra crossed module. 
Section $3$ gives the definition of principal $2$-bundles with which we work.
It is taken from Wockel's article \cite{Woc}, together with restrictions from
\cite{BaeSch}.
In Section $4$, we discuss in general $L_{\infty}$-valued differential forms 
on the manifold $M$, based on the article of Getzler \cite{Get}. We believe 
that this is the right generalization of the calculus of Lie algebra valued
differential forms needed for ordinary principal $G$-bundles. We find a 
curious $3$-form term (see equation (\ref{**})) 
in the Maurer-Cartan equation for differential forms 
with values in a semi-strict Lie $2$-algebra which seems to be new. 
In Section $5$, we construct the connection $1$-form ${\mathcal A}_0$
of Baez-Schreiber from the band of the non-abelian gerbe. It is not so clear 
in \cite{BaeSch} on which differential geometric object the construction of
${\mathcal A}_0$ is carried out, and we believe that expressing it as the usual
iterated integral construction on the band (which is an ordinary principal 
$G$-bundle !) is of conceptual importance. 

Section $6$ is the heart of our article and explains the mechanism to transform
the flat connection ${\mathcal A}_0$ into a Hochschild cycle for the 
differential graded algebra $CH_*(\Omega^*,\Omega^*)$. It lives therefore
in the Hochschild homology of the algebra of Hochschild chains. 
Section $7$ recalls from \cite{GTZ} that ``Hochschild of Hochschild''-homology
is the higher Hochschild homology associated to the torus $\T^2$.
   
The main theorem of the present article is the construction of the 
homology cycle:

\begin{theo}   \label{main_theorem}
Consider a non-abelian principal $2$-bundle with trivial band 
on a simply connected manifold $M$ 
with a structure crossed module
$\mu:{\mathfrak h}\to{\mathfrak g}$ such that the Lie algebra
${\mathfrak h}$ is abelian. Then the connection $1$-form ${\mathcal A}_0$
of Baez-Schreiber gives rise to a cycle $P({\mathcal A}_0)$ in the higher 
Hochschild homology $HH_{\bullet}^{\T}(\Omega^*,\Omega^*)$ which corresponds 
to the holonomy of the gerbe.
\end{theo}   

As stated before, we do not consider the condition that ${\mathfrak h}$ is 
abelian as a restriction of generality, because up to equivalence,
it may be achieved for an arbitrary crossed module.
 
By construction, the cycle $P({\mathcal A}_0)$ is not always trivial, i.e. 
a boundary. The triviality condition on the band may be understood as 
expressing that the construction is local. The gluing of the locally defined
connection $1$-forms of Baez and Schreiber to a global connection $1$-form
(see \cite{BaeSch}) should permit to glue our Hochschild cycles 
$P({\mathcal A}_0)$ to a global cycle. This will be taken up in futur work.

Another subject of further research is to understand that the connection
$1$-form ${\mathcal A}_0$ does not only lead to a higher Hochschild cycle
w.r.t. the $2$-dimensional torus, but actually to higher Hochschild cycles
w.r.t. any compact topological surface. In fact, we believe that
there is a way to recover
$HH_{\bullet}^{\Sigma_g}$ for a connected compact surface $\Sigma_g$ of 
genus $g$ from $HH_{\bullet}^{\T}$. 

\vspace{.5cm}     

\noindent{\bf Acknowledgements:}
We thank Gregory Ginot for answering a question about higher Hochschild 
homology. We also thank Urs 
Schreiber and Dimitry Roytenberg for answering a question in the context
of Section $1$. 

\section{Strict Lie $2$-algebras and crossed modules}
\label{Section_one}

We gather in this section preliminaries on strict Lie $2$-algebras and 
crossed modules, and their relation to semistrict Lie $2$-algebras. 
The main result is the possibility to replace a crossed module 
$\mu:{\mathfrak h}\to{\mathfrak g}$ by an equivalent one having abelian
${\mathfrak h}$. This will be important for defining holonomy as a cycle in
higher Hochschild homology.  

Lie $2$-algebras have been
the object of different studies, see \cite{BaeCra} for semi-strict Lie 
$2$-algebras or \cite{Roy} for (general weak) Lie $2$-algebras. 

\subsection{Strict $2$-vector spaces}

We fix a field $\K$ of characteristic $0$; in geometrical situations, we will 
always take $\K=\R$. A {\it $2$-vector space} $V$ over $\K$ is simply a 
category object in ${\tt Vect}$, the category of vector spaces (cf Def. $5$ in 
\cite{BaeCra}). This means that $V$ consists of a vector space of {\it arrows} 
$V_{-1}$, a vector space of {\it objects} $V_0$, linear maps 
$\xymatrix{V_{-1}\ar@<2pt>[r]^{s}\ar@<-2pt>[r]_{t} & V_0}$ called {\it source}
and {\it target}, a linear map $i:V_0\to V_{-1}$, called {\it object 
inclusion}, and a linear map 
$$m\,:\,V_{-1}\times_{V_0}V_{-1}\to V_{-1},$$
which is called the {\it categorical composition}.
These data is supposed to satisfy the usual axioms of a category.

An equivalent point of view is to regard a $2$-vector space as
a $2$-term complex of vector spaces $d:C_{-1}\to C_0$. Pay attention
to the change in degree with respect to \cite{BaeCra}. We use here a 
cohomological convention, instead of their homological convention, in order
to have the right degrees for the differential forms with values in crossed 
modules later on.   

The equivalence between $2$-vector spaces and $2$-term complexes is spelt 
out in Section $3$ of \cite{BaeCra}:
one passes from a category object in ${\tt Vect}$ (given by
$\xymatrix{V_{-1}\ar@<2pt>[r]^{s}\ar@<-2pt>[r]_{t} & V_0}$, 
$i:V_0\to V_{-1}$ etc) to a $2$-term complex $d:C_{-1}\to C_0$ by taking
$C_{-1}:=\ker(s)$, $d:=t|_{\ker(s)}$ and $C_0=V_0$. In the reverse direction,
to a given $2$-term complex $d:C_{-1}\to C_0$, one associates 
$V_{-1}=C_0\oplus C_{-1}$, $V_0=C_0$, $s(c_0,c_{-1})=c_0$, 
$t(c_0,c_{-1})=c_0+d(c_{-1})$, and $i(c_0)=(c_0,0)$. The only subtle point
is here that the categorical composition $m$ is already determined by 
$\xymatrix{V_{-1}\ar@<2pt>[r]^{s}\ar@<-2pt>[r]_{t} & V_0}$ and 
$i:V_0\to V_{-1}$ (see Lemma $6$ in \cite{BaeCra}). Namely, writing an
arrow $c_{-1}=:f$ with $s(f)=x$, $t(f)=y$, i.e. $f:x\mapsto y$, one denotes
the {\it arrow part} of $f$ by $\vec{f}:=f-i(s(f))$, and for two 
composable arrows
$f,g\in V_{-1}$, the composition $m$ is then defined by
$$f\circ g\,:=\,m(f,g)\,:=\,i(x)+\vec{f}+\vec{g}.$$ 

\subsection{Strict Lie $2$-algebras and crossed modules} 

\begin{defi}
A strict Lie $2$-algebra is a category object in the category ${\tt Lie}$ of 
Lie algebras over $\K$. 
\end{defi}

This means that it is the data of two Lie algebras, 
${\mathfrak g}_0$, the {\it Lie algebra of objects}, and ${\mathfrak g}_{-1}$, 
the {\it Lie algebra of arrows}, together with morphisms of Lie algebras
$s,t:{\mathfrak g}_{-1}\to {\mathfrak g}_0$, {\it source} and {\it target}, 
a morphism $i:{\mathfrak g}_0\to {\mathfrak g}_{-1}$, the 
{\it object inclusion},
and a morphism $m:{\mathfrak g}_{-1}\times_{{\mathfrak g}_0}{\mathfrak g}_{-1}
\to {\mathfrak g}_{-1}$, the {\it composition of arrows}, such that the usual 
axioms of a category are satisfied. 

Let us now come to crossed modules of Lie algebras. We refer to \cite{Wag} for
more details. 

\begin{defi}
A crossed module of Lie algebras is a morphism of Lie algebras
$\mu:{\mathfrak h}\to{\mathfrak g}$ together with an action of ${\mathfrak g}$ 
on ${\mathfrak h}$ by derivations such that for all $h,h'\in{\mathfrak h}$
and all $g\in{\mathfrak g}$
\begin{itemize}
\item[(a)] $\mu(g\cdot h)\,=\,[g,\mu(h)]$ and
\item[(b)] $\mu(h)\cdot h'\,=\,[h,h']$.
\end{itemize}
\end{defi}

One may associate to a crossed module of Lie algebras a $4$-term exact 
sequence of Lie algebras
$$0\to V\to{\mathfrak h}\stackrel{\mu}{\to}{\mathfrak g}\to\bar{\mathfrak g}
\to 0,$$
where we used the notation $V:=\ker(\mu)$ and $\bar{\mathfrak g}:=\coker(\mu)$.
It follows from the properties (a) and (b) of a crossed module that
$\mu({\mathfrak h})$ is an ideal, so $\bar{\mathfrak g}$ is a Lie algebra, and 
that $V$ is a central ideal of ${\mathfrak h}$ and a 
$\bar{\mathfrak g}$-module (because the outer action, to be 
defined below, is a genuine action on the center of ${\mathfrak h}$).   

Recall the definition of the {\it outer action}
$s:\bar{\mathfrak g}\to{\mathfrak o}{\mathfrak u}{\mathfrak t}({\mathfrak h})$
for a crossed module of Lie algebras $\mu:{\mathfrak h}\to{\mathfrak g}$. 
The Lie algebra 
$${\mathfrak o}{\mathfrak u}{\mathfrak t}({\mathfrak h}):=
{\mathfrak d}{\mathfrak e}{\mathfrak r}({\mathfrak h})/{\rm ad}({\mathfrak h})
$$
is the Lie algebra of outer derivations of ${\mathfrak h}$, i.e. the quotient 
of the Lie algebra of all derivations 
${\mathfrak d}{\mathfrak e}{\mathfrak r}({\mathfrak h})$ by the ideal
${\rm ad}({\mathfrak h})$ of inner derivations, i.e. those of the form
$h'\mapsto[h,h']$ for some $h\in{\mathfrak h}$. 
  
To define $s$, choose a
linear section $\rho:\bar{\mathfrak g}\to{\mathfrak g}$ and compute its default
to be a homomorphism of Lie algebras
$$\alpha(x,y):=[\rho(x),\rho(y)]-\rho([x,y]),$$
for $x,y\in\bar{\mathfrak g}$. As the projection onto $\bar{\mathfrak g}$ is 
a homomorphism of Lie algebras, $\alpha(x,y)$ is in its kernel, and there 
exists therefore an element $\beta(x,y)\in{\mathfrak h}$ such that
$\mu(\beta(x,y))=\alpha(x,y)$. 

We have for all $h\in{\mathfrak h}$
$$\big(\rho(x)\circ\rho(y)-\rho(y)\circ\rho(x)-\rho([x,y])\big)\cdot h=
 \alpha(x,y)\cdot h=\mu(\beta(x,y))\cdot h=[\beta(x,y),h],$$
and in this sense, elements of $\bar{\mathfrak g}$ act on ${\mathfrak h}$ up
to inner derivations. We obtain a well defined homomorphism of Lie algebras
$$s:\bar{\mathfrak g}\to{\mathfrak o}{\mathfrak u}{\mathfrak t}({\mathfrak h})
$$
by $s(x)(h)=\rho(x)\cdot h$.  
 
Strict Lie $2$-algebras are in one-to-one correspondance
with crossed modules of Lie algebras, like in the case of groups, cf 
\cite{Loday}. For the convenience of the reader, let us include this here:

\begin{theo} \label{theorem_crmod_strict_Lie}
Strict Lie $2$-algebras are in one-to-one correspondence
with crossed modules of Lie algebras.
\end{theo}

\pr Given a Lie $2$-algebra 
$\xymatrix{{\mathfrak g}_{-1}\ar@<2pt>[r]^{s}\ar@<-2pt>[r]_{t} & 
{\mathfrak g}_0}$, $i:{\mathfrak g}_0\to{\mathfrak g}_{-1}$,
the corresponding crossed module is defined by
$$\mu:=t|_{\ker(s)}\,:\,{\mathfrak h}:=\ker(s)\to{\mathfrak g}:=
{\mathfrak g}_0.$$
The action of ${\mathfrak g}$ on ${\mathfrak h}$ is given by
$$g\cdot h\,:=\,[i(g),h],$$
for $g\in{\mathfrak g}$ and $h\in{\mathfrak h}$ (where the bracket is taken in 
${\mathfrak g}_{-1}$). This is well defined and an action by derivations.
Axiom (a) follows from 
$$\mu(g\cdot h)=\mu([i(g),h])=[\mu\circ i(g),\mu(h)]=[g,\mu(h)].$$
Axiom (b) follows from
$$\mu(h)\cdot h'=[i\circ\mu(h),h']=[i\circ t(h),h']$$
by writing $i\circ t(h)=h+r$ for $r\in\ker(t)$ and by using that $\ker(t)$ and 
$\ker(s)$ in a Lie $2$-algebra commute (shown in Lemma 
\ref{lemma_commuting_kernels} after the proof).

On the other hand, given a crossed module of Lie algebras 
$\mu:{\mathfrak h}\to{\mathfrak g}$, associate to it
$$\xymatrix{{\mathfrak h}\rtimes{\mathfrak g}\ar@<2pt>[r]^{s}\ar@<-2pt>[r]_{t} 
& {\mathfrak g}},\,\,\,\,i:{\mathfrak g}\to{\mathfrak h}\rtimes{\mathfrak g}$$
by $s(h,g)=g$, $t(h,g)=\mu(h)+g$, $i(g)=(0,g)$, where the semi-direct product
Lie algebra ${\mathfrak h}\rtimes{\mathfrak g}$ is built from the given action 
of ${\mathfrak g}$ on ${\mathfrak h}$. Let us emphasize that 
${\mathfrak h}\rtimes{\mathfrak g}$ is built from the Lie algebra 
${\mathfrak g}$ and the ${\mathfrak g}$-module ${\mathfrak h}$; the bracket of
${\mathfrak h}$ does not intervene here. The composition of arrows is already
encoded in the underlying structure of $2$-vector space, as remarked in the 
previous subsection. \fin

\begin{lem} \label{lemma_commuting_kernels}
$[\ker(s),\ker(t)]\,=\,0$ in a Lie $2$-algebra.
\end{lem}    

\pr The fact that the composition of arrows is a homomorphism of Lie algebras
gives the following ``middle four exchange'' (or functoriality) property
$$[g_1,g_2]\circ[f_1,f_2]\,=\,[g_1\circ f_1,g_2\circ f_2]$$
for composable arrows $f_1,f_2,g_1,g_2\in{\mathfrak g}_1$. 
Now suppose that $g_1\in\ker(s)$ and $f_2\in\ker(t)$. Then denote by $f_1$ and 
by $g_2$ the
identity (w.r.t. the composition) in $0\in{\mathfrak g}_0$. As these are 
identities, we have $g_1=g_1\circ f_1$ and $f_2=g_2\circ f_2$. On the other 
hand, $i$ is a morphism of Lie algebras and sends $0\in{\mathfrak g}_0$ to
the $0\in{\mathfrak g}_1$. Therefore we may conclude
$$[g_1,f_2]\,=\,[g_1\circ f_1,g_2\circ f_2]\,=\,[g_1,g_2]\circ[f_1,f_2]\,=\,
0.$$\fin  

Furthermore, it is well-known (cf \cite{Wag}) that (equivalence classes of) 
crossed modules of Lie algebras are classified by third cohomology classes.

\begin{rem}
It is implicit in the previous proof that starting from a crossed module
$\mu:{\mathfrak h}\to{\mathfrak g}$, passing to the Lie $2$-algebra 
$\xymatrix{{\mathfrak g}_{-1}\ar@<2pt>[r]^{s}\ar@<-2pt>[r]_{t} & 
{\mathfrak g}_0}$, $i:{\mathfrak g}_0\to{\mathfrak g}_{-1}$ (and thus 
forgetting
the bracket on ${\mathfrak h}$ !), one may finally reconstruct the bracket 
on ${\mathfrak h}$. This is due to the fact that it is encoded in the action 
and the morphism, using the property (b) of a crossed module.
\end{rem}

\subsection{Semi-strict Lie $2$-algebras and 
$2$-term $L_{\infty}$-algebras}

An equivalent point of view is to regard a strict 
Lie $2$-algebra as a Lie algebra object in the category ${\tt Cat}$ of (small) 
categories. From this second point of view, we have a functorial Lie bracket 
which is supposed to be antisymmetric and must fulfill the Jacobi identity.
Weakening the antisymmetry axiom and the Jacobi identity up to coherent 
isomorphisms leads then to semi-strict Lie $2$-algebras (here antisymmetry
holds strictly, but Jacobi is weakened), hemi-strict Lie $2$-algebras 
(here Jacobi holds strictly, but antisymmetry is weakened) or even to (general)
Lie $2$-algebras (both axioms are weakened). Let us record the definition 
of a semi-strict Lie $2$-algebra (see \cite{BaeCra} Def. $22$):

\begin{defi}
A semi-strict Lie $2$-algebra consists a $2$-vector space $L$ together with
a skew-symmetric, bilinear and functorial bracket $[,]:L\times L\to L$ and
a completely antisymmetric trilinear natural isomorphism 
$$J_{x,y,z}\,:\,[[x,y],z]\to[x,[y,z]]+[[x,z],y],$$
called the {\it Jacobiator}. The Jacibiator is required to satisfy the 
Jacobiator identity (see \cite{BaeCra} Def. $22$).
\end{defi}
  
Semi-strict Lie $2$-algebras together with morphisms of semi-strict Lie 
$2$-algebras (see Def. $23$ in \cite{BaeCra}) form a strict $2$-category
(see Prop. $25$ in \cite{BaeCra}). 
Strict Lie algebras form a full sub-$2$-category of this $2$-category, see
Prop. $42$ in \cite{BaeCra}. In order to regard a strict Lie $2$-algebra
${\mathfrak g}_{-1}\to {\mathfrak g}_0$ as a semi-strict Lie $2$-algebra, the 
functorial bracket is constructed for $f:x\mapsto y$ and $g:a\mapsto b$, 
$f,g\in
{\mathfrak g}_{-1}$ and $x,y,a,b\in{\mathfrak g}_0$ by defining its source
$s([f,g])$ and its arrow part $\vec{[f,g]}$ to be $s([f,g]):=[x,a]$ and
$\vec{[f,g]}:=[x,\vec{g}]+[\vec{f},b]$ (see proof of Thm. $36$ in 
\cite{BaeCra}). By construction, it is compatible with the composition, i.e.
functorial. 

\begin{rem}
One observes that the functorial bracket on a strict Lie $2$-algebra
${\mathfrak g}_{-1}\to {\mathfrak g}_0$ is constructed from the bracket in
${\mathfrak g}_0$, and the bracket between ${\mathfrak g}_{-1}$ and 
${\mathfrak g}_0$, but does not involve the bracket on  ${\mathfrak g}_{-1}$
itself. 
\end{rem}

There is a $2$-vector space underlying every semi-strict Lie $2$-algebras, thus
one may ask which structure is inherited from a semi-strict Lie $2$-algebra
by the corresponding $2$-term complex of vector spaces. This leads us to
$2$-term $L_{\infty}$-algebras, see \cite{BaeCra} Thm. $36$. Our definition
here differs from theirs as we stick to the cohomological setting and degree
$+1$ differentials, see \cite{Get} Def. $4.1$. 

\begin{defi}
An $L_{\infty}$-algebra is a graded vector space $L$ together with a sequence 
$l_k(x_1,\ldots,x_k)$, $k>0$, of graded antisymmetric operations of degree 
$2-k$ such that the following identity is satisfied:

\begin{equation*}   
 \sum_{k=1}^n(-1)^k\sum_{\begin{array}{c} 
\scriptstyle i_1< \ldots < i_k;\,j_1 < \ldots < j_{n-k} \\ 
\scriptstyle \{i_1,\ldots,i_k\}\cup\{j_1,\ldots,j_{n-k}\}=\{1,\ldots,n\}
\end{array}}
(-1)^{\epsilon}l_n(l_k(x_{i_1},\ldots,x_{i_k}),x_{j_1},\ldots,x_{j_{n-k}})=0.
\end{equation*}
Here, the sign $(-1)^{\epsilon}$ equals the product of the sign of the shuffle 
permutation and the Koszul sign. 
\end{defi}

We will be mainly concerned with $2$-term $L_{\infty}$-algebras. These are 
$L_{\infty}$-algebras $L$ such that the graded vector space $L$ consists only
of two components $L_0$ and $L_{-1}$. An $L_{\infty}$-algebra 
$L=L_0\oplus L_{-1}$ has at most $l_1$, $l_2$ and $l_3$ as its non-trivial 
``brackets''. $l_1$ is a differential (i.e. here just a linear map 
$L_0\to L_{-1}$), $l_2$ is a bracket with components $[,]:L_0\otimes L_0\to 
L_0$ and $[,]:L_{-1}\otimes L_0\to L_{-1}$, $[,]:L_0\times L_{-1}\to L_{-1}$, 
and $l_3$ is some kind of $3$-cocycle $l_3:L_0\otimes L_0\otimes L_0\to 
L_{-1}$. More precisely, in case $l_1=0$, $L_0$ is a Lie algebra, $L_{-1}$
is an $L_0$-module and $l_3$ is then an actual $3$-cocycle. This kind of
$2$-term $L_{\infty}$-algebra is called {\it skeletal}, see Section $6$
in \cite{BaeCra} and our next subsection.  The complete 
axioms satisfied  by $l_1$, $l_2$ and $l_3$ in a $2$-term $L_{\infty}$-algebra
are listed in Lemma $33$ of \cite{BaeCra}. 

As said before, the passage from a $2$-vector space to its associated $2$-term
complex induces a passage from semi-strict Lie $2$-algebras to
$2$-term $L_{\infty}$-algebras, which turns out to be an equivalence of $2$ 
categories (see \cite{BaeCra} Thm. $36$):

\begin{theo} 
The $2$-categories of semi-strict Lie $2$-algebras and of 
$2$-term $L_{\infty}$-algebras are equivalent.
\end{theo}

\begin{rem}
In particular, restricting to the sub-$2$-category of strict Lie 
$2$ algebras, there is an equivalence between crossed modules of Lie algebras 
and $2$-term $L_{\infty}$-algebras with trivial $l_3$. In other words, there 
is an equivalence between crossed modules and differential graded Lie algebras.
\end{rem}

\subsection{Classification of semi-strict Lie $2$-algebras}

Baez and Crans show in \cite{BaeCra} that every semi-strict Lie $2$-algebra
is equivalent to a skeletal Lie $2$-algebra (i.e. one where the differential 
$d$ of the underlying complex of vector spaces is zero). Then they go on by
showing that skeletal Lie $2$-algebras are classified by triplets consisting 
of an honest Lie algebra $\bar{\mathfrak g}$, a $\bar{\mathfrak g}$-module 
$V$ and a
class $[\gamma]\in H^3(\bar{\mathfrak g},V)$. 
This is achieved using the homotopy
equivalence of the underlying complex of vector spaces with its cohomology.
In total, they get in this way a classification, up to equivalence, of 
semi-strict Lie $2$-algebras in terms of triplets 
$(\bar{\mathfrak g},V,[\gamma])$.

On the other hand, strict Lie $2$-algebras are in one-to-one correspondance
with crossed modules of Lie algebras, as we have seen in a previous 
subsection.
In conclusion, there are two ways to classify strict  Lie $2$-algebras: by
the associated crossed module or, regarding them as special semi-strict Lie 
$2$-algebras, by Baez-Crans classification. 
Let us show here that these two classifications are compatible, i.e. that
they lead to the same triplet $(\bar{\mathfrak g},V,[\gamma])$. 

For this, let us 
denote by ${\tt sLie2}$ the $2$-category of strict Lie $2$-algebras, by 
${\tt ssLie2}$ the $2$-category of semi-strict Lie $2$-algebras, by 
${\tt sssLie2}$ the $2$-category of skeletal semi-strict Lie $2$-algebras, by
${\tt triplets}$ the (trivial) $2$-category of triplets of the above form 
$({\mathfrak g},V,[\gamma])$, and by ${\tt crmod}$ the $2$-category 
of crossed modules of Lie algebras. 

\begin{theo}
The following diagram is commutative:
     
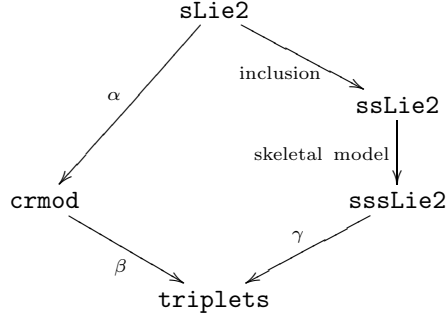
\vspace{.3cm}
\hspace{2cm}
\xymatrix{
            & {\tt sLie2} \ar[dr]_{\rm inclusion}\ar[ddl]_{\alpha}& \\
            &             & {\tt ssLie2} \ar[d]_{\rm skeletal\,\,\, model}
\\
{\tt crmod} \ar[dr]_{\beta}&             & {\tt sssLie2} \ar[dl]_{\gamma}
\\
            &{\tt triplets}& }

The $2$-functors $\alpha$ and $\gamma$ are bijections, while the $2$-functor 
$\beta$ induces a bijection when passing to equivalence classes. 
\end{theo}

\pr Let us first describe the arrows. The arrow 
$\alpha:{\tt sLie2}\to{\tt crmod}$
has been investigated in Theorem \ref{theorem_crmod_strict_Lie}. The arrow 
$\beta:{\tt crmod}\to{\tt triplets}$ sends a crossed module 
$\mu:{\mathfrak h}\to{\mathfrak g}$ to the triple 
$$(\coker(\mu)=:\bar{\mathfrak g},\ker(\mu)=:V,[\gamma]),$$
where the cohomology class $[\gamma]\in H^3(\bar{\mathfrak g},V)$
is defined choosing sections -- the procedure is described in detail in 
\cite{Wag}. The arrow ${\tt ssLie2}\to{\tt sssLie2}$ is the choice of a
skeletal model for a given semi-strict Lie $2$-algebra -- it is given by the
homotopy equivalence of the underlying $2$-term complex with its cohomology
displayed in the extremal lines of the following diagram

\vspace{.3cm}
\hspace{3.5cm}
\xymatrix{ C_{-1} \ar[r]^d & C_0 \ar@{=}[d] \\
\ker(d) \ar@{=}[d] \ar@{^{(}->}[u] \ar[r]^0 & C_0 \ar@{->>}[d] \\
\ker(d) \ar[r]^0 & C_0\,/\,\im(d)}
\vspace{.3cm}

The arrow $\gamma: {\tt sssLie2}\to{\tt triplets}$ sends a skeletal $2$-Lie 
algebra to the triplet defined by the cohomology class of $l_3$ (cf 
\cite{BaeCra}). 

Now let us show that the diagram commutes. For this,
let $d:C_{-1}\to C_0$ with some bracket $[,]$ and $l_3=0$ be a
$2$-term $L_{\infty}$-algebra corresponding to seeing a strict Lie 
$2$-algebra as a semi-strict Lie $2$-algebra, and build its skeletal model.
The model comes together with a morphism of semi-strict Lie $2$-algebras 
$(\phi_{2},\phi_{-1},\phi_0)$ given by 

\vspace{.3cm}
\hspace{3.5cm}
\xymatrix{ C_{-1} \ar[r]^d & C_0  \\
\ker(d)  \ar@{^{(}->}[u]^{\phi_{-1}} \ar[r]^0 & C_0\,/\,\im(d)  
\ar[u]^{\phi_0}  }
\vspace{.3cm}  

Here $\phi_0=:\sigma$ is a linear section of the quotient map. 
The structure of a
semi-strict Lie $2$-algebra is transfered to the lower line in order to make  
$(\phi_{2},\phi_{-1},\phi_0)$ a morphism of semi-strict Lie $2$-algebras.
In order to compute now the $l_3$ term of the lower semi-strict Lie 
$2$-algebra, one first finds
that (first equation in definition 34 of \cite{BaeCra})
$\phi_{2}:C_0\,/\,\im(d)\times C_0\,/\,\im(d)\to C_{-1}$ is such that
$$d\phi_{2}(x,y)\,=\,\sigma[x,y] - [\sigma(x),\sigma(y)],$$
the default of the section $\sigma$ to be a homomorphism of Lie algebras.
Then $l_3$ is related to $\phi_{2}$ by the second formula in definition 34 of
\cite{BaeCra}. This gives here
$$l_3(x,y,z)\,=\,(d_{\rm CE}\phi_{2})(x,y,z)$$
for $x,y,z\in C_0\,/\,\im(d)$. $d_{\rm CE}$ is the formal Chevalley-Eilenberg 
differential of the cochain 
$\phi_{2}:C_0\,/\,\im(d)\times C_0\,/\,\im(d)
\to C_{-1}$ with values in $C_{-1}$ as if $C_{-1}$ was a 
$C_0\,/\,\im(d)$-module (which is usually not the case). This is exactly the
expression of the cocycle $\gamma$ associated to the crossed module of Lie 
algebras $d:C_{-1}\to C_0$ obtained using the section $\sigma$, see 
\cite{Wag}.\fin

\begin{cor}  \label{cor_equivalence_to_strict}
Every semi-strict Lie $2$-algebra is equivalent (as an object of the 
$2$-category ${\tt ssLie2}$) to a strict Lie $2$-algebra.
\end{cor}

This corollary is already known because of abstract reasons. 
Here we have 
proved a result somewhat more refined: the procedure to strictify a 
semi-strict Lie $2$-algebra is rather easy to perform. First one has to
pass to cohomology by homotopy equivalence, and then one has to construct
the crossed module corresponding to a given cohomology class. This can be done 
in several way, using free Lie algebras \cite{LodKas}, using injective modules
\cite{Wag} etc. and one may adapt the construction method to the problem at 
hand.

\subsection{The construction of an abelian representative} 

We will show in this section that to a given class 
$[\gamma]\in H^3(\bar{\mathfrak g},V)$,
there exists a crossed module of Lie algebras 
$\mu:{\mathfrak h}\to{\mathfrak g}$ with class $[\gamma]$ (and $\ker(\mu)=V$ 
and $\coker(\mu)=\bar{\mathfrak g}$) such that 
${\mathfrak h}$ is abelian. This will be important for the treatment in higher
Hochschild homology of the holonomy of a gerbe.   
 
\begin{theo}  \label{theorem_reduction_to_abelian}
For any $[\gamma]\in H^3(\bar{\mathfrak g},V)$, there exists a crossed 
module of Lie algebras 
$\mu:{\mathfrak h}\to{\mathfrak g}$ with associated class $[\gamma]$
such that $\ker(\mu)=V$, $\coker(\mu)=\bar{\mathfrak g}$ and ${\mathfrak h}$ 
is abelian.
\end{theo}

\pr This is Theorem 3 in \cite{Wag}. Let us sketch its proof here. The category
of $\bar{\mathfrak g}$-modules has enough injectives, therefore $V$ may be
embedded in an injective $\bar{\mathfrak g}$-module $I$. 
We obtain a short exact 
sequence of $\bar{\mathfrak g}$-modules
$$0\to V\stackrel{i}{\to} I\stackrel{\pi}{\to} Q\to 0,$$
where $Q:=I/V$ is the quotient. $I$ injective implies 
$H^p(\bar{\mathfrak g},I)=0$
for all $p>0$. Therefore the short exact sequence of coefficients induces 
a connective homomorphism
$$\partial:H^2(\bar{\mathfrak g},Q)\to H^3(\bar{\mathfrak g},V)$$
which is an isomorphism. To $[\gamma]$ corresponds thus a class
$[\alpha]\in H^2(\bar{\mathfrak g},Q)$ with $\partial[\alpha]=[\gamma]$. A 
representative $\alpha\in Z^2(\bar{\mathfrak g},Q)$ gives rise to an abelian 
extension
$$0\to Q\to Q\times_{\alpha}\bar{\mathfrak g}\to\bar{\mathfrak g}\to 0.$$
Now one easily verifies (see the proof of Theorem 3 in \cite{Wag}) that 
the splicing together of the short exact coefficient sequence and the abelian 
extension gives rise to a crossed module
$$0\to V\to I\to Q\times_{\alpha}\bar{\mathfrak g}\to\bar{\mathfrak g}\to 0.$$
More precisely, the crossed module is 
$\mu:I\to Q\times_{\alpha}\bar{\mathfrak g}$ given by $\mu(x)=(\pi(x),0)$, 
the action of ${\mathfrak g}:=Q\times_{\alpha}\bar{\mathfrak g}$ on
${\mathfrak h}:=I$ is induced by the action of $\bar{\mathfrak g}$ on $I$ 
and the Lie bracket is trivial on $I$, i.e. $I$ is abelian. 

One also easily verifies (see the proof of Theorem 3 in \cite{Wag}) 
that the associated cohomology class for such a crossed module (which is 
the Yoneda product of a short exact coefficient sequence and an abelian 
extension) is $\partial[\alpha]$, the image under the connective homomorphism
(induced by the short exact coefficent sequence) of the class defining the 
abelian extension. Therefore the associated class is here 
$\partial[\alpha]=[\gamma]$ as required.
\fin

We thus obtain the following refinement of Corollary 
\ref{cor_equivalence_to_strict}:  

\begin{cor}  \label{cor_reduction}
Every semi-strict Lie $2$-algebra is equivalent to a strict Lie $2$-algebra 
corresponding to a crossed module $\mu:{\mathfrak h}\to{\mathfrak g}$ with
abelian ${\mathfrak h}$, such that ${\mathfrak h}$ is a 
$\bar{\mathfrak g}:={\mathfrak g}/\mu({\mathfrak h})$-module and such that
the outer action is a genuine action.
\end{cor}

\pr This follows from Corollary \ref{cor_equivalence_to_strict} together with
Theorem \ref{theorem_reduction_to_abelian}. The fact that ${\mathfrak h}$ is a 
$\bar{\mathfrak g}=:{\mathfrak g}/\mu({\mathfrak h})$-module and that the 
outer action is a genuine action are equivalent. They are true either
by inspection of the representative constructed in the proof of Theorem 
\ref{theorem_reduction_to_abelian}, or by the following argument:

The outer action $s$ is an action only up to inner derivations. 
But these are trivial in case ${\mathfrak h}$ is abelian:
$$\mu(h)\cdot h'=[h,h']=0$$
for all $h,h'\in{\mathfrak h}$ by property (b) of a crossed module. 
\fin 
 
\begin{rem}
An analoguous statement is true on the level of (abstract) groups and even 
topological groups \cite{WagWoc}.  
Unfortunately, we ignore whether such
a statement is true in the category of Lie groups, i.e. given 
a locally smooth group $3$-cocycle $\gamma$ on $\bar{G}$ with values in a 
smooth $\bar{G}$-module $V$, is there a smooth (split) crossed 
module of Lie groups $\mu:H\to G$ with $H$ abelian and cohomology class
$[\gamma]$ ? From the point of view of Lie algebras, there are two steps 
involved: having solved the problem on the level of Lie algebras (as above),
one has to integrate the $2$-cocycle $\alpha$. This is well-understood
thanks to work of Neeb. The (possible) obstructions lie in $\pi_1(\bar{G})$ and 
$\pi_2(\bar{G})$, and vanish thus for simply connected, finite dimensional 
Lie groups $\bar{G}$. 
The second step is to integrate the involved $\bar{\mathfrak g}$-module
$I$ to a $\bar{G}$-module. 
As $I$ is necessarily infinite dimensional, this is the 
hard part of the problem.
\end{rem}

\section{Crossed modules of Lie groups}

In this section, we introduce the strict Lie $2$-groups which will be the 
typical fiber of our principal $2$-bundles. While the notion of a crossed 
module of groups is well-understood and purely algebraic, the notion of a 
crossed module of Lie groups involves subtle smoothness requirements.  

We will heavily draw on \cite{Neeb} and adopt Neeb's point of view, namely, we 
regard a crossed module of Lie groups as a central extension $\hat{N}\to N$ of 
a normal split Lie subgroup $N$ in a Lie group $G$ for which the conjugation 
action of $G$ on $N$ lifts to a smooth action on $\hat{N}$. This point of
view is linked to the one regarding a crossed module as a homomorphism
$\mu:H\to G$ by taking $H=\hat{N}$ and $\im(\mu)=N$. 

\begin{defi}
A morphism of Lie groups $\mu:H\to G$, together with a homomorphism 
$\hat{S}:G\to{\rm Aut}(H)$ defining a smooth action $\hat{S}:G\times H\to
H$, $(g,h)\mapsto g\cdot h=\hat{S}(g)(h)$ of $G$ on $H$, is called a (split)
crossed module of Lie groups if the following conditions are satisfied:
\begin{enumerate}
\item $\mu\circ\hat{S}(g)={\rm conj}_{\mu(g)}\circ\mu$ for all $g\in G$.
\item $\hat{S}\circ\mu:H\to{\rm Aut}(H)$ is the conjugation action.
\item $\ker(\mu)$ is a split Lie subgroup of $H$ and $\im(\mu)$ is a split
Lie subgroup of $G$ for which $\mu$ induces an isomorphism $H/\ker(\mu)\to
\im(\mu)$. 
\end{enumerate}
\end{defi} 

Recall that in a split crossed module of Lie groups $\mu:H\to G$, the quotient 
Lie group $\bar{G}:=G/\mu(H)$ acts smoothly (up to inner automorphisms) on $H$. 
This outer action $S$ of $\bar{G}$ on $H$ is a homomorphism 
$S:\bar{G}\to{\rm Out}(H)$ which is constructed like in the case of Lie 
algebras. The smoothness of $S$ follows directly from the splitting 
assumptions. Here ${\rm Out}(H)$ denotes the 
{\it group of outer automorphisms} of $H$, defined by
$${\rm Out}(H):={\rm Aut}(H)/{\rm Inn}(H),$$
where ${\rm Inn}(H)\subset{\rm Aut}(H)$ is the normal subgroup of automorphisms
of the form $h'\mapsto hh'h^{-1}$ for some $h\in H$.  

It is shown in \cite{Neeb} that one may associate to a (split)
crossed module of Lie groups a locally smooth $3$-cocycle $\gamma$
(whose class is the obstruction against the realization of the outer action 
in terms of an extension).
 
It is clear that a (split) crossed module of Lie groups induces a crossed 
module of the corresponding Lie algebras.

\begin{defi} \label{equivalent_crossed_modules}
Two crossed modules $\mu:M\to N$ (with action $\eta$) and $\mu':M'\to N'$
(with action $\eta'$) such that $\ker(\mu)=\ker(\mu')=:V$ and 
$\coker(\mu)=\coker(\mu')=:G$ are called {\it elementary equivalent} if there
are group homomorphisms $\varphi:M\to M'$ and $\psi:N\to N'$ which are 
compatible with the actions, i.e.
$$\varphi(\eta(n)(m))=\eta'(\psi(n))(\varphi(m))$$
for all $n\in N$ and all $m\in M$, and such that the following diagram is
commutative: 

\hspace{2cm}
\xymatrix{
0 \ar[r] & V \ar[d]^{\id_V} \ar[r]^{i} & M \ar[d]^{\varphi} 
\ar[r]^{\mu} & N \ar[d]^{\psi} \ar[r]^{\pi} &  G 
\ar[d]^{\id_G} \ar[r] & 0 \\
0 \ar[r] & V  \ar[r]^{i'} & M' \ar[r]^{\mu'} & N'  
\ar[r]^{\pi'} &  G  \ar[r] & 0}
\vspace{.5cm}

We call {\it equivalence of crossed modules} the equivalence relation generated
by elementary equivalence. One easily sees that two crossed modules are 
equivalent in case there exists a zig-zag of elementary equivalences going
from one to the other (where the arrows do not necessarily all go into the same
direction). 
\end{defi}

In the context of crossed modules of Lie groups, all morphisms are supposed 
to be morphisms of Lie groups, i.e. smooth.

\section{Principal $2$-bundles and gerbes}

In this section, we will start introducing the basic geometric objects of our
study, namely principal $2$-bundles and gerbes. 
We choose to work here with a strict Lie $2$-group ${\mathcal G}$, i.e. a 
split crossed module of Lie groups, and its 
associated crossed module of Lie algebras $\mu:{\mathfrak h}\to
{\mathfrak g}$, and to consider principal $2$-bundles and gerbes which are 
defined by non-abelian cocycles (or transition functions). 
The principal object which we will use later on is the band of a gerbe.

\subsection{Definition}

In order to keep notations and abstraction to a reasonable minimum, we will 
consider geometric objects like bundles, gerbes, etc only over an honest 
(finite dimensional) base
manifold $M$, instead of considering a ringed topos, a stack or anything else.
                 
Let $\mu:H\to G$ be a (split) crossed module of Lie groups. Let our base space
$M$ be an honest (ordinary) manifold, and let ${\mathcal U}=\{U_i\}$ be a 
good open cover of $M$.   
The following definition is based on p.29 in \cite{BaeSch}
and on the corresponding presentation in \cite{Woc}). 

\begin{defi}
A non-abelian cocycle $(g_{ij},h_{ijk})$ is the data of (smooth) transition 
functions
$$g_{ij}:U_i\cap U_j\to G$$
and 
$$h_{ijk}:U_i\cap U_j\cap U_k\to H$$
which satisfy the non-abelian cocycle identities
$$\mu(h_{ijk}(x))g_{ij}(x)g_{jk}(x)=g_{ik}(x)$$
for all $x\in U_{ijk}:=U_i\cap U_j\cap U_k$, and 
$$h_{ikl}(x)h_{ijk}(x)=h_{ijl}(x)(g_{ij}(x)\cdot h_{jkl}(x))$$
for all $x\in  U_{ijkl}:=U_i\cap U_j\cap U_k\cap U_l$.
\end{defi}

The \v{C}ech cochains $g_{ij}$ and $h_{ijk}$ are (by definition) antisymmetric 
in the indices.
One may complete the set of indices to all pairs resp. triplets by imposing
the functions to be equal to $1_G$ resp $1_H$ on repeated indices. 

We go on by defining equivalence of non-abelian cocycles with values in the 
same crossed module of Lie groups $\mu:H\to G$: 

\begin{defi}
Two non-abelian cocycles $(g_{ij},h_{ijk})$ and $(g_{ij}',h_{ijk}')$ on the 
same cover are said to be {\it equivalent} if there exist (smooth) functions
$\gamma_i:U_i\to G$ and $\eta_{ij}:U_{ij}\to H$ such that
$$\gamma_i(x)g_{ij}'(x)=\mu(\eta_{ij}(x))g_{ij}(x)\gamma_j(x)$$
for all $x\in U_{ij}$, and 
$$\eta_{ik}(x)h_{ijk}(x)=(\gamma_i(x)\cdot h_{ijk}'(x))\eta_{ij}(x)(g_{ij}(x)
\cdot \eta_{jk}(x))$$
for all $x\in U_{ijk}$.
\end{defi}
In general, one should define equivalence for cocycles 
corresponding to different covers. Passing to a common refinement, one easily 
adapts the above definition to this framework (this is spelt out in
\cite{Woc}). 

\begin{defi}
A principal $2$-bundle, also called (non-abelian) gerbe and denoted 
${\mathcal G}$, is the data of an equivalence class of non-abelian cocycles.
\end{defi}

By abuse of language, we will also call a representative $(g_{ij},h_{ijk})$
a principal $2$-bundle or a (non-abelian) gerbe.    

\begin{lem}
If the (split) crossed module of Lie groups 
$\mu:H\to G$ is
replaced by an equivalent crossed module $\mu':H'\to G'$, then a given 
non-abelian cocycle $(g_{ij},h_{ijk})$ taking values in $\mu:H\to G$ gives rise
to a non-abelian cocycle $(g_{ij}',h_{ijk}')$ taking values in 
$\mu':H'\to G'$.
\end{lem}

\pr This is rather formal. We restrict here to an elementary equivalence - 
for an arbitrary equivalence, one should iterate the argument.

Given a non-abelian cocycle $(g_{ij},h_{ijk})$ and an elementary equivalence 
$(\phi,\psi):(H,G)\to(H',G')$, define a non-abelian cocycle 
$(g_{ij}',h_{ijk}')$ by $g_{ij}':=\psi(g_{ij}')$ and  $h_{ijk}':=\phi(h_{ijk}')$.
It is easily checked that $(g_{ij}',h_{ijk}')$ satisfies the non-abelian
cocycle conditions thanks to the requirement that $(\phi,\psi):(H,G)\to(H',G')$
is an elementary equivalence.\fin 

\begin{rem}
Let us denote by $(\phi,\psi)_*(g_{ij},h_{ijk})$ the thus constructed cocycle
$(g_{ij}',h_{ijk}')$ (w.r.t. the elementary equivalence $(\phi,\psi)$).
It would make sense to define that a non-abelian cocycle $(g_{ij},h_{ijk})$ 
w.r.t. the crossed module $\mu:H\to G$ is {\it elementary equivalent} to 
a non-abelian cocycle $(g_{ij}',h_{ijk}')$ w.r.t. the (possibly different
but elementary equivalent) 
crossed module $\mu':H'\to G'$ in case the cocycle $(g_{ij}',h_{ijk}')$ 
is equivalent to $(\phi,\psi)_*(g_{ij},h_{ijk})$ 
as cocycles w.r.t. $\mu':H'\to G'$ (where $(\phi,\psi)$
is the elementary equivalence from $\mu:H\to G$ to $\mu':H'\to G'$). 
One may then use this relation of elementary equivalence to define
(arbitrary) equivalence between cocycles.
\end{rem}  

Recall that for a split crossed module of Lie groups $\mu:H\to G$, the 
image $\mu(H)$ is a normal Lie subgroup of $G$, and the quotient group
$\bar{G}:=G/\mu(H)$ is therefore a Lie group. 

\begin{lem}
Let ${\mathcal G}$ be a gerbe defined by the cocycle $(g_{ij},h_{ijk})$.

Then one may associate to ${\mathcal G}$ an ordinary
principal $\bar{G}$-bundle ${\mathcal B}$ on $M$ which has as its 
transition functions the composition of the $g_{ij}$ and the canonical 
projection $G\to G/\mu(H)=\bar{G}$.
\end{lem}

\pr This is clear. Indeed, passing to the quotient $G\to G/\mu(H)$, 
the identity
$$\mu(h_{ijk}(x))g_{ij}(x)g_{jk}(x)=g_{ik}(x)$$  
becomes the cocycle identity
$$\bar{g}_{ij}(x)\bar{g}_{jk}(x)=\bar{g}_{ik}(x)$$
for a principal $G/\mu(H)$-bundle on $M$ defined by the transition functions
$$\bar{g}_{ij}:U_{ij}\to G/\mu(H)$$
obtained from composing $g_{ij}:U_{ij}\to G$ with the projection 
$G\to G/\mu(H)$.\fin
   
\begin{defi}
The principal $\bar{G}$-bundle ${\mathcal B}$ on $M$ associated to 
the gerbe ${\mathcal G}$ defined by the cocycle $(g_{ij},h_{ijk})$ is
called the band of ${\mathcal G}$. 
\end{defi}

\subsection{Connection data}  \label{section_connection_data}

Let, as before, $M$ be a manifold and let ${\mathcal U}=\{U_i\}$ be a 
good open cover of $M$. Let ${\mathcal G}$ be a gerbe defined by the 
cocycle $(g_{ij},h_{ijk})$. We associate to ${\mathcal G}$ now connection data
like in \cite{BaeSch} Sect. 2.1.4. 

\begin{defi}
{\it Connection data} for the non-abelian cocycle $(g_{ij},h_{ijk})$ is the data
of {\it connection $1$-forms} $A_i\in\Omega^1(U_i,{\mathfrak g})$ and of 
{\it curving $2$-forms} $B_i\in\Omega^2(U_i,{\mathfrak h})$, together with
{\it connection transformation $1$-forms} $a_{ij}\in\Omega(U_{ij},{\mathfrak h})$
and {\it curving transformation $2$-forms} 
$d_{ij}\in\Omega^2(U_{ij},{\mathfrak h})$ such that the following laws hold:
\begin{enumerate}
\item[(a)] transition law for connection $1$-forms on $U_{ij}$
$$A_i+\mu(a_{ij})=g_{ij}A_jg_{ij}^{-1}+g_{ij}dg_{ij}^{-1}.$$
\item[(b)] transition law for the curving $2$-forms on $U_{ij}$
$$B_i=g_{ij}\cdot B_j+da_{ij}.$$
\item[(c)] transition law for the curving transformation $2$-forms on $U_{ijk}$
$$d_{ij}+g_{ij}\cdot d_{jk}=h_{ijk}d_{ik}h_{ijk}^{-1}+
h_{ijk}(\mu(B_i)+F_{A_i})h_{ijk}^{-1}.$$
\item[(d)] coherence law for the transformers of connection $1$-forms on 
$U_{ijk}$
$$0=a_{ij}+g_{ij}\cdot a_{jk}-h_{ijk} a_{ik}h_{ijk}^{-1}-h_{ijk}dh_{ijk}^{-1}-h_{ijk}
(A_i\cdot h_{ijk}^{-1}).$$
\end{enumerate}
\end{defi}
 
In accordance with \cite{BaeSch} equation (2.73) on p. 59, we will choose 
$d_{ij}=0$ in the following. The transition law $(c)$ for the curving 
transformation $2$-forms reads then simply   
$$0=\mu(B_i)+F_{A_i},$$
which is the equation of {\it vanishing fake curvature}. In the following,
we will always suppose that the fake curvature vanishes (cf Section 
\ref{L_infinity_valued_forms}).  

\begin{defi}  \label{definition_curvature_3_form}
Given a ${\mathcal G}$ be a gerbe defined by the 
cocycle $(g_{ij},h_{ijk})$ with connection data $(A_i,B_i,a_{ij})$, the 
{\it curvature $3$-form} $H_i\in\Omega^3(U_i,{\mathfrak h})$ is defined by
$$H_i=d_{A_i}B_i,$$
i.e. it is the covariant derivative of the curving $2$-form
$B_i\in\Omega^2(U_i,{\mathfrak h})$
with respect to the connection $1$-form $A_i\in\Omega^1(U_i,{\mathfrak h})$.

Its transformation law on $U_{ij}$ is  
$$H_i=g_{ij}\cdot H_j,$$
(because in our setting fake curvature and curving transformation $2$-forms 
vanish). 
\end{defi}

Observe that only the crossed module of Lie algebras $\mu:{\mathfrak h}\to
{\mathfrak g}$ plays a role as values of the differential forms.
According to Section \ref{Section_one}, it constitutes no restriction of 
generality (up to equivalence) to consider ${\mathfrak h}$ abelian. 
In our main application (construction of the holonomy higher Hochschild 
cycle), we will suppose ${\mathfrak h}$ to be abelian. Many steps on the way 
are true for arbitrary ${\mathfrak h}$. 
The property of being abelian simplifies the above coherence law for the 
transformers of connection 
$1$-forms on $U_{ijk}$ for which we thus obtain in the abelian setting:
$$0=a_{ij}+g_{ij}\cdot a_{jk}-a_{ik}-h_{ijk}dh_{ijk}^{-1}-h_{ijk}
(A_i\cdot h_{ijk}^{-1}).$$

We note in passing:
\begin{lem}
The connection $1$-forms induce an ordinary connection on the band 
${\mathcal B}$ of the gerbe ${\mathcal G}$.
\end{lem}

\pr This follows at once from equation $(a)$ in the definition of connection 
data.\fin

On the other hand, we will always be in a local setting, therefore, in the 
following, we will drop the indices $i,j,k,\ldots$ which refer to the open set
we are on. 

\section{$L_{\infty}$-valued differential forms}
\label{L_infinity_valued_forms}

In this section, we will associate to each principal $2$-bundle an 
$L_{\infty}$-algebra of $L_{\infty}$-valued differential forms. 
This $L_{\infty}$-algebra
replaces the differential graded Lie algebra of Lie algebra valued forms which 
plays a role for ordinary principal $G$-bundles. Here, the  
$L_{\infty}$-algebra of values (of the differential forms) will be the 
$2$-term $L_{\infty}$-algebra associated to the strict structure Lie 
$2$-algebra of the principal $2$-bundle.   
We follow closely \cite{Get}, section $4$.

Given an $L_{\infty}$-algebra ${\mathfrak g}_{\infty}$ and a manifold $M$, 
the tensor 
product $\Omega^*(M)\otimes{\mathfrak g}_{\infty}$ of 
${\mathfrak g}_{\infty}$-valued smooth 
differential forms on $M$ is an $L_{\infty}$-algebra by prolongating the 
$L_{\infty}$-operations of ${\mathfrak g}_{\infty}$ point by point to differential 
forms. The only point to notice is that the de Rham differential 
$d_{\rm de Rham}$ gives a 
contribution to the first bracket 
$l_1:{\mathfrak g}_{\infty}\to{\mathfrak g}_{\infty}$ which is 
also a differential of degree $1$. 

We will apply this scheme to the $2$-term $L_{\infty}$-algebras arising from
a semi-strict Lie $2$-algebra 
${\mathfrak g}_{\infty}=({\mathfrak g}_{-1},{\mathfrak g}_{0})$.
The only (possibly) non-zero operations are the 
differential $l_1$, the bracket $[,]=l_2$ and the $3$-cocycle $l_3$.
Our choice of degrees is that an element 
$\alpha^k\in\Omega^*(M)\otimes{\mathfrak g}_{\infty}$ is of degree $k$ in case
$\alpha^k\in\bigoplus_{i\geq 0}\Omega^i(M)\otimes {\mathfrak g}_{k-i}$.
An element of degree $1$ is thus a sum $\alpha^1=\alpha_1+\alpha_2$ with
$\alpha_1\in\Omega^1(M)\otimes {\mathfrak g}_{0}$ and 
$\alpha_2\in\Omega^2(M)\otimes {\mathfrak g}_{-1}$. 

Recall the following definitions (cf Def. $4.2$ in \cite{Get}):

\begin{defi}
The Maurer--Cartan set ${\rm MC}({\mathfrak g}_{\infty})$ of a nilpotent 
$L_{\infty}$-algebra ${\mathfrak g}_{\infty}$ is the set of 
$\alpha\in{\mathfrak g}^1$
satisfying the Maurer--Cartan equation ${\mathcal F}(\alpha)=0$. More 
explicitely, this means 
$${\mathcal F}(\alpha)\,:=\,l_1\alpha+\sum_{k=2}^{\infty}\frac{1}{k!}l_k(
\alpha,\ldots,\alpha)\,=\,0.$$
\end{defi}

The Maurer--Cartan 
equations for the degree $1$ elements of the $L_{\infty}$-algebra 
$\Omega^*(M)\otimes{\mathfrak g}_{\infty}$ (see \cite{Get} Def. 
$4.3$) read therefore
\begin{equation} \label{*}
d_{\rm de Rham}\alpha_1+\frac{1}{2}[\alpha_1,\alpha_1]+l_1\alpha_2\,=\,0,
\end{equation}
and
\begin{equation}  \label{**}
d_{\rm de Rham}\alpha_2+[\alpha_1,\alpha_2] + 
l_3(\alpha_1,\alpha_1,\alpha_1)\,=\,0.
\end{equation}
Equation (\ref{*}) is an equation of $2$-forms; in the gerbe literature
it is known as the equation of the vanishing of the fake curvature.
Equation (\ref{**}) is an equation of $3$-forms and seems to be new in this 
context. The special case $l_3=0$ corresponds to example $6.5.1.3$ in 
\cite{SSS}, when
one interpretes $d_{\rm de Rham}\alpha_2+[\alpha_1,\alpha_2]$ as
the covariant derivative $d_{\alpha_1}\alpha_2$. When applied to connection
data of a non-abelian gerbe (see Section \ref{section_connection_data}),
the vanishing of the covariant derivative means that the $3$-curvature
(cf Definition \ref{definition_curvature_3_form}) of the gerbe vanishes. 
This is 
sometimes expressed as being a {\it flat gerbe}.   

Let us record the special case of a strict Lie $2$-algebra 
${\mathfrak g}_{\infty}$ given by a crossed module 
$\mu:{\mathfrak h}\to{\mathfrak g}$ for later use:

\begin{lem} \label{Maurer-Cartan_equation}
A degree $1$ element of the $L_{\infty}$-algebra 
$\Omega^*(M)\otimes{\mathfrak g}_{\infty}$ is a pair $(A,B)$ with 
$A\in\Omega^1(M)\otimes{\mathfrak g}$ and 
$B\in\Omega^2(M)\otimes{\mathfrak h}$.

The element $(A,B)$ satisfies the Maurer-Cartan equation if and only if 
$$d_{\rm de Rham}A+\frac{1}{2}[A,A]+\mu B\,=\,0\,\,\,\,{\rm and}\,\,\,\, 
d_{\rm de Rham}B+[A,B]=0.$$
\end{lem}

Elements of degree $0$ in $\Omega^*(M)\otimes{\mathfrak g}$ are 
sums $\alpha^0=\beta_0+\beta_1$ with 
$\beta_0\in\Omega^0(M)\otimes {\mathfrak g}_{0}$ and 
$\beta_1\in\Omega^1(M)\otimes {\mathfrak g}_{-1}$. These act by gauge 
transformations on elements of the Maurer-Cartan set. Namely, $\beta_0$
has to be exponentiated to an element $B_0\in\Omega^0(M,G_0)$ (where $G_0$ is 
the connected, $1$-connected Lie group corresponding to ${\mathfrak g}_{0}$)
and leads then 
to gauge transformations of the first kind, in the sense of \cite{BaeSch}.
Elements $\beta_1$ lead directly to gauge transformations of the second kind, 
in the sense of \cite{BaeSch}. The fact that they don't have to be 
exponentiated corresponds to the fact that there is no bracket on 
the ${\mathfrak g}_{-1}$-part of the $L_{\infty}$-algebra.
These gauge transformations will not play a role in the present paper, but
will become a central subject when gluing the local expressions of the 
connection $1$-form of Baez-Schreiber to a global connection.  

\begin{defi}
Let ${\mathfrak g}$ be a nilpotent $L_{\infty}$-algebra. 
The Maurer--Cartan variety ${\mathcal M}{\mathcal C}({\mathfrak g})$ 
is the quotient of the 
Maurer--Cartan set ${\rm MC}({\mathfrak g})$ by the exponentiated action
of the infinitesimal automorphisms ${\mathfrak g}^0$ of 
${\rm MC}({\mathfrak g})$.
\end{defi}

We do not assert that the quotient ${\mathcal M}{\mathcal C}({\mathfrak g})$
is indeed a variety. It is considered here as a set.

\section{Path space and the connection $1$-form associated to 
a principal  $2$-bundle}

In this section we explain how the connective structure on a 
gerbe gives rise to a connection on path space.

\subsection{Path space as a Fr\'echet manifold}

We first recall some basic facts about path spaces which allow us
to employ the basic notion of differential geometry, in particular 
differential forms and connections. For a manifold 
$M$, let $\cP M:=C^\infty([0,1],M )$  be the space of paths in $M$.
Baez and Schreiber \cite{BaeSch} fix in their definition 
the starting point and the end point of the paths, i.e.
for two points $s$ and $t$ in $M$, $\cP_s^t M$ denotes the space
of paths from $s$ (``source'') to $t$ (``target''). 
The space $\cP M$ can made into a Fr\'echet manifold modeled
on the Fr\'echet space $C^\infty([0,1],\R^n)$ in case $n$ is the dimension 
of $M$. 
Similar constructions exhibit the loop space  
$$LM:=\{\gamma:[0,1]\stackrel{{\mathcal C}^{\infty}}{\rightarrow} M\,|\,\, 
\gamma(0)=\gamma(1),\,\,\,\gamma^{(k)}(0)=\gamma^{(k)}(1)=0\,\,\,\forall 
k\geq 1\}$$ 
as a Fr\'echet manifold, 
see for example \cite{Nee} for a detailed account on the Fr\'echet 
manifold structure on (this version of) $L M$ in case $M$ is a Lie group.
The generalization to arbirary $M$ is quite standard. Let us emphasize
that this version of $LM$ comes to mind naturally when writing the circle 
$S^1$ as $[0,1]/\sim$ in ${\mathcal C}^{\infty}(S^1,M)$. The fact that one
demands $\gamma^{(k)}(0)=\gamma^{(k)}(1)=0$ and not only 
$\gamma^{(k)}(0)=\gamma^{(k)}(1)$ for all $k\geq 1$ is sometimes expressed by 
saying that the loops have a ``sitting instant''.      
 
For differentiable 
Fr\'echet manifolds, one can introduce differential forms, de Rham 
differential and prove a De Rham theorem for smoothly paracompact Fr\'echet 
manifolds. The only thing beyond the necessary definitions that we need from 
Fr\'echet differential geometry is an expression of the de Rham differential on
$LM$, an expression due to Chen \cite{Chen} which will play its role in 
the proof of Proposition \ref{prop_curvature}.

\subsection{The connection $1$-form of Baez-Schreiber}

Let $\mu:H\to G$ be a split crossed module of Lie groups. Denote by $S$ 
the outer action of $\bar{G}$ on $H$, i.e. the
homomorphism $S:\bar{G}\to{\rm Out}(H)$.  

Composing the transition functions $\bar{g}_{ij}:U_{ij}\to \bar{G}$ with the 
homomorphism $S:\bar{G}\to{\rm Out}(H)$, we obtain the transition functions
of an ${\rm Out}(H)$-principal bundle denoted ${\mathcal B}_S$. This is
then an ordinary principal bundle, and we may apply ordinary holonomy
theory to the principal bundle ${\mathcal B}_S$. 

As we are only interested in these 
constructions and these constructions are purely on Lie algebra level, we 
will neglect now the crossed module of Lie groups and focus on the crossed 
module of Lie algebras. In doing so, we may assume (up to equivalence without 
loss of generality) that in the crossed module
$\mu:{\mathfrak h}\to{\mathfrak h}$, ${\mathfrak h}$ is abelian and that
the outer action 
$s:\bar{\mathfrak g}\to{\mathfrak o}{\mathfrak u}{\mathfrak t}({\mathfrak h})$
(associated to $\mu:{\mathfrak h}\to{\mathfrak g}$ like in Section 1) 
is a genuine action (see Corollary \ref{cor_reduction}).

In the following, we will suppose that the principal bundle ${\mathcal B}_S$
is trivial (or, in other words, we will do a local construction).
A connection $1$-form on ${\mathcal B}_S$ is then simply 
a differential form $A_S\in\Omega^1_{\End}$, and given 
a $1$-form $A\in\Omega^1(M,\bar{\mathfrak g})$, one obtains such a form
$A_S$ by $A_S:=s\circ A$. We will suppose that ${\mathcal B}_S$ posseses
a flat connection $\nabla$ which will be our reference point in the affine 
space of connections.   

\begin{rem}
Actually, in case the $1$-form $A\in\Omega^1(M,{\mathfrak g})$ (and not in 
$\Omega^1(M,\bar{\mathfrak g})$ !), there is no problem 
to define the action of $A$ on $B$. We do not need ${\mathfrak h}$ to be 
abelian here (in case we do not want to use the band, for example).
\end{rem}

Consider now the loop space $LM$ of $M$.  
Let us procede with the Wilson loop or iterated integral 
construction of Section 6 of 
\cite{AbbZei}. For every $n\geq 0$, consider the $n$-simplex
$$\triangle^n:=\{(t_0,t_1,\ldots,t_n,t_{n+1})\,|\,0=t_0\leq t_1\leq\ldots\leq
t_n\leq t_{n+1}=1\}.$$
Define the evaluation maps ${\rm ev}$ and ${\rm ev}_{n,i}$ as follows:
\begin{eqnarray*}
&&{\rm ev}:\triangle^n\times LM\to M \\
&&{\rm ev}(t_0,t_1,\ldots,t_n,t_{n+1};\gamma)=\gamma(0)=\gamma(1)\\
&&{\rm ev}_{n,i}:\triangle^n\times LM\to M \\
&&{\rm ev}_{n,i}(t_0,t_1,\ldots,t_n,t_{n+1};\gamma)=\gamma(t_i).
\end{eqnarray*}
Denote by ${\rm ad}\,{\mathcal B}_s$ the adjoint bundle associated
to the principal bundle ${\mathcal B}_S$ using the adjoint action of 
${\rm Out}(H)$ on ${\mathfrak o}{\mathfrak u}{\mathfrak t}({\mathfrak h})
={\mathfrak d}{\mathfrak e}{\mathfrak r}({\mathfrak h})$.
Let $T_i:{\rm ev}_{n,i}^*({\rm ad}\,{\mathcal B}_s)\to{\rm ev}^*({\rm ad}\,
{\mathcal B}_s)$ denote the map, between pullbacks of adjoint bundles to
$\triangle^n\times LM$, defined at a point $(0=t_0,t_1,\ldots,t_n,t_{n+1}=1;
\gamma)$ by the parallel transport along and in the direction of $\gamma$ 
from $\gamma(t_i)$ to $\gamma(t_{n+1})=\gamma(1)$ in the bundle 
${\rm ad}\,{\mathcal B}_s$ with respect to the flat connection $\nabla$.

For $\alpha_i\in \Omega^*(M,{\rm ad}\,{\mathcal B}_s)$, $1\leq i\leq n$, define
$\alpha(n,i)\in \Omega^*(\triangle^n\times LM,{\rm ev}^*{\rm ad}\,
{\mathcal B}_s)$ by
$$\alpha(n,i)=T_i{\rm ev}_{n,i}^*\alpha_i.$$
The associated bundle to ${\mathcal B}_s$ with typical fiber the universal 
enveloping algebra $U{\mathfrak d}{\mathfrak e}{\mathfrak r}({\mathfrak h})$
is denoted by ${\mathcal B}_s^U$. Now define $V_{\alpha_1,\ldots,\alpha_n}^n\in
\Omega^*(LM,{\rm ev}^*{\mathcal B}_s^U)$ by
\begin{eqnarray*}
&V^0=1,& \\
&V_{\alpha_1,\ldots,\alpha_n}^n=\int_{\triangle^n}\alpha(n,1)\wedge\ldots
\wedge\alpha(n,n)\,\,\,{\rm for}\,\,\,n\geq 1,&
\end{eqnarray*}
and set
$$V_{\alpha}=\sum_{n=0}^{\infty}V_{\alpha}^n,\,\,\,{\rm where}\,\,\,
V_{\alpha}^n=V_{\alpha,\ldots,\alpha}^n.$$
It is noteworthy that this infinite sum is convergent. This is shown in 
\cite{AbbZei} in Appendix B. Observe that for $1$-forms 
$\alpha_1,\ldots,\alpha_n$, the loop space form $V_{\alpha_1,\ldots,\alpha_n}^n$ 
has degree $0$ for all $n$. 

Furthermore, define for $B\in\Omega^2(M,{\mathfrak h})$ and $\sigma\in[0,1]$
the $1$-form $B^*(\sigma)\in\Omega^1(LM,{\mathfrak h})$ by 
$$B^*(\sigma)=i_K{\rm EV}^*_{\sigma}B$$
for the evaluation map ${\rm EV}_{\sigma}:LM\to M$, $EV(\gamma):=
\gamma(\sigma)$ and the 
vector field $K$ on $LM$ which is the infinitesimal generator of the 
$S^1$-action on $LM$ by rigid rotations.  

Now fix an element $(A,B)$ of the Maurer-Cartan set w.r.t. some Lie
algebra crossed module $\mu:{\mathfrak h}\to{\mathfrak g}$. 
Evaluating elements of $\End({\mathfrak h})$ on ${\mathfrak h}$, we obtain 
a connection $1$-form ${\mathcal A}_0$ on $LM$ with values in 
${\mathfrak h}$ given by
$${\mathcal A}_0=\int_0^1V_{A}(B^*(\sigma))d\sigma.$$
(Indeed, as $A$ is a $1$-form, the loop space form $V_A$ is of degree $0$, and 
$V_{A}(B^*(\sigma))$ is of degree $1$ and remains of degree $1$ after 
integration w.r.t. $\sigma$.)  

This gives the formula for the connection $1$-form of Baez and Schreiber
on p. 43 of \cite{BaeSch}:

\begin{prop}
The constructed connection $1$-form ${\mathcal A}_0$ on $LM$ with values in 
${\mathfrak h}$ coincides with the path space $1$-form ${\mathcal A}_{(A,B)}=
\oint_A(B)$ of Def. 2.23 in \cite{BaeSch}.
\end{prop}

\pr This follows from a step-by-step comparison.\fin

\section{The holonomy cycle associated to a principal $2$-bundle}

A central construction of \cite{ATZ} associates to elements ${\mathcal A}$ in 
the Maurer-Cartan space a holonomy class $[P({\mathcal A})]$ in 
$HH_*(\Omega^*,\Omega^*)$.
This is done using the following proposition (cf {\it loc. cit.} Section 4):

\begin{prop} \label{prop_Hochschild_cycle}
Suppose given a differential graded associative algebra $\Omega^*$ and
an element ${\mathcal A}\in\Omega^{\rm odd}$.  
The following are equivalent.
\begin{itemize}
\item[(a)] ${\mathcal A}$ is a Maurer-Cartan element, i.e. $d{\mathcal A}+
{\mathcal A}\cdot {\mathcal A}=0$,
\item[(b)] the chain 
$$P({\mathcal A})\,:=\,1\otimes 1 + 1\otimes {\mathcal A} + 1\otimes 
{\mathcal A}\otimes {\mathcal A}+\ldots$$
in the Hochschild complex $CH_*(\Omega^*,\Omega^*)$ is a cycle. 
\end{itemize}
\end{prop} 
 
\pr We have
$$d_{\rm Hoch}(P({\mathcal A}))=\sum \pm 1\otimes{\mathcal A}\ldots
\otimes d{\mathcal A}\ldots\otimes{\mathcal A}+\sum 
\pm 1\otimes{\mathcal A}\ldots\otimes {\mathcal A}\cdot{\mathcal A}
\ldots\otimes{\mathcal A}.$$
Therefore the cycle property is equivalent to 
$d{\mathcal A}+{\mathcal A}\cdot{\mathcal A}=0$, i.e. to the Maurer-Cartan 
equation.\fin

The degrees are taken such that all terms in $P({\mathcal A})$ are of degree
$0$ in case ${\mathcal A}$ is of degree $1$, i.e. the degrees of $\Omega^*$
are shifted by one. This is the correct degree when taking Hochschild homology 
as a model for loop space cohomology. 
 
We will apply this proposition to the connection $1$-form  ${\mathcal A}_0$ on 
$LM$. The $1$-form ${\mathcal A}_0$ is an element of 
$\Omega^1(LM,U{\mathfrak h})$. 
The condition that ${\mathcal A}_0$ is a Maurer-Cartan element is then
that the curvature of ${\mathcal A}_0$ vanishes. This curvature has been 
computed in \cite{BaeSch} p. 43 to be given by the following formula
(needless to say, no assumption is made on ${\mathfrak h}$ for this 
computation):

\begin{prop} \label{prop_curvature}
The curvature of the $1$-form ${\mathcal A}_0$ is equal to 
\begin{eqnarray*}
{\mathcal F}_{{\mathcal A}_0}&:=&-\oint_A(d_AB)-\oint_A(d\alpha(T_a)(B),
(F_A+\mu(B))^a)\\
&:=&
\int_0^1V_A((d_AB)^*(\sigma))d\sigma-\int_0^1V_A((F_A+\mu(B))^*(\sigma))d\sigma,
\end{eqnarray*}
where $d_AB$ is the covariant derivative of $B$ w.r.t. $A$ and 
$F_A+\mu(B)$ is the fake curvature of the couple $(A,B)$.
\end{prop}

\pr A detailed proof is given in \cite{BaeSch} Prop. 2.7 and Cor. 2.2, 
p. 42-43. Here we will
only sketch the main steps of the proof.

First compute $d_{\rm deRham}{\mathcal A}_0$ for the de Rham differential 
$d_{\rm deRham}$. As explained
in {\it loc. cit.} Prop. 2.4, p. 35, the action of the de Rham differential 
on a Chen form $\oint_A(\omega_1,\ldots,\omega_n)$ is given by two terms, 
namely 
$$\sum_k\pm\oint_A(\omega_1,\ldots,d_{\rm deRham}\omega_k,\ldots,\omega_n)$$ 
and 
$$\sum_k\pm\oint_A(\omega_1,\ldots,\omega_{k-1}\wedge\omega_k,\ldots,\omega_n)
.$$
In our case, we get thus four terms, according to whether the $B$ is 
involved or not. The terms which does not involve $B$ give a term involving
the curvature $F_A=dA+A\cdot A$. The terms involving $B$ give a term
involving the covariant derivative $d_AB=d_{\rm deRham}B+A\cdot B$ of $B$ 
w.r.t. $A$.

Now the computation of the curvature of ${\mathcal A}_0$ adds to the 
de Rham derivative $d_{\rm deRham}{\mathcal A}_0$ a term 
${\mathcal A}_0\cdot{\mathcal A}_0$. This term is easily seen to 
be the term involving $\mu(B)$.\fin    
 
We see that for a Maurer-Cartan element $(A,B)$ (in the sense of Section
\ref{L_infinity_valued_forms}), ${\mathcal A}_0$ is a flat 
connection (by Lemma \ref{Maurer-Cartan_equation} and Proposition 
\ref{prop_curvature}), and 
$P({\mathcal A}_0)$ is therefore a Hochschild cycle 
(by Proposition \ref{prop_Hochschild_cycle}). 

\begin{defi}
The Hochschild cycle $P({\mathcal A}_0)$ is the
holonomy cycle associated to the given principal $2$-bundle.
\end{defi}  

Let us abbreviate $\Omega^*(M,U{\mathfrak h})$ to $\Omega^*$. 
Our main point is now that the assumption that ${\mathfrak h}$ is abelian
implies that $\Omega^*$ (and for the same reason also 
$\Omega^*(LM,U{\mathfrak h})$) is a 
{\it commutative} differential graded algebra, thus the shuffle product endows
the (ordinary) Hochschild complex
$CH_*(\Omega^*,\Omega^*)$ with the structure
of a differential graded commutative algebra (cf \cite{Loday} Cor. 4.2.7, 
p. 125). On the other hand, we have for a 
simply connected manifold $M$:

\begin{lem}
There is an quasi-isomorphism of commutative differential graded algebras 
$$\Omega^*(LM,U{\mathfrak h})\simeq CH_*(\Omega^*,\Omega^*).$$
\end{lem}

\pr Let us first observe that the loop space with sitting instant is 
homotopically equivalent to $LM={\mathcal C}^{\infty}(S^1,M)$. Indeed, the
version with sitting instant corresponds to taking $[0,1]/\sim$ as $S^1$. 
Elements for a rigorous proof (at least in the case of Lie groups $M$)
may be found in Appendix A of \cite{Nee}. 
 
Therefore, our assertion is a version with coefficients in the graded 
associative algebra 
$U{\mathfrak h}$ of Corollary 2.6, p. 11, in \cite{Lod2}, originally shown by 
Chen \cite{Chen}. Observe that the 
coefficients do not contribute to the differentials.
\fin 

In conclusion, we obtain a homology class 
$$[P({\mathcal A}_0)]\in HH_*(CH_*(\Omega^*,\Omega^*),CH_*(\Omega^*,\Omega^*)).$$
In the next section, we will explain how to interprete 
$$HH_*(CH_*(\Omega^*,\Omega^*),CH_*(\Omega^*,\Omega^*))$$ 
in terms of higher Hochschild homology
as $HH_*^{\T}(\Omega^*,\Omega^*)$, the {\it higher Hochschild homology}
of the $2$-dimensional torus $\T$. We therefore obtain
$$[P({\mathcal A}_0)]\in HH_*^{\T}(\Omega^*,\Omega^*).$$

\section{Higher Hochschild homology}

In this section, we consider higher Hochschild homology. It has been 
invented by Pirashvili in \cite{Pir} and further developed by Ginot, Tradler 
and Zeinalian in \cite{GTZ}. Here, we follow closely \cite{GTZ}.

In order to define higher Hochschild homology, it is essential to 
restrict to {\it commutative} differential graded associative algebras 
$\Omega^*$. We will see below explicitely why this is the case.

Denote by $\triangle$ the (standard) category whose objects are the finite
ordered sets
$[k]=\{0,1,\ldots,k\}$ and morphisms $f:[k]\to[l]$ are non-decreasing maps,
i.e. for $i>j$, one has $f(i)\geq f(j)$. 
Special non-decreasing maps are the injections $\delta_i:[k-1]\to[k]$ 
characterized by missing $i$ (for $i=0,\ldots,k$)
and the surjections $\sigma_j:[k]\to[k-1]$ which send $j$ and $j+1$ to $j$ 
(equally for $j=0,\ldots,k$).

A {\it finite simplicial set}
$Y_{\bullet}$ is by definition a contravariant functor $Y_{\bullet}:
\triangle^{\rm op}\to{\tt Sets}$. 
The sets of $k$-simplices are denoted $Y_k:=Y([k])$. 
The induced maps $d_i:=Y_{\bullet}(\delta_i)$ and $s_j:=Y_{\bullet}(\sigma_j)$
are called {\it faces} and {\it degeneracies} respectively. 
Let $Y_{\bullet}$ be a pointed finite simplicial set. 
For $k\geq 0$, we put $y_k:=|Y_k|-1$, i.e. 
one less than the cardinal of the finite set $Y_k$.

The {\it higher Hochschild chain complex} of $\Omega^*$ associated to the 
simplicial set $Y_{\bullet}$ (and with values in $\Omega^*$) is 
defined by 
$$CH_{\bullet}^{Y_{\bullet}}(\Omega^*,\Omega^*):=\bigoplus_{n\in\Z}
CH_n^{Y_{\bullet}}(\Omega^*,\Omega^*),$$
where
$$CH_n^{Y_{\bullet}}(\Omega^*,\Omega^*):=\bigoplus_{k\geq 0}(\Omega^*\otimes
(\Omega^*)^{\otimes y_k})_{n+k}.$$
In order to define the differential, define induced maps as follows. 
For any map $f:Y_k\to Y_l$ of pointed sets
and any (homogeneous) element $m\otimes a_1\otimes\ldots\otimes a_{y_k}\in
\Omega^*\otimes(\Omega^*)^{\otimes y_k}$, we denote by $f_*:
\Omega^*\otimes(\Omega^*)^{\otimes y_k}\to\Omega^*\otimes(\Omega^*)^{\otimes y_l}$
the map
$$f_*(m\otimes a_1\otimes\ldots\otimes a_{y_k}):=(-1)^{\epsilon}
n\otimes b_1\otimes\ldots\otimes b_{y_l},$$
where $b_j=\Pi_{i\in f^{-1}(j)}a_i$ (or $b_j=1$ in case $f^{-1}(j)=\emptyset$)
for $j=0,\ldots,y_l$, and $n=m\cdot\Pi_{i\in f^{-1}({\rm basepoint}), 
j\not={\rm basepoint}}a_i$. The sign $\epsilon$ is determined by the usual 
Koszul sign rule. The above face and degeneracy maps $d_i$ and $s_j$ induce
thus boundary maps $(d_i)_*:CH_k^{Y_{\bullet}}(\Omega^*,\Omega^*)\to
CH_{k-1}=k^{Y_{\bullet}}(\Omega^*,\Omega^*)$ and degeneracy maps 
$(s_j)_*:CH_{k-1}^{Y_{\bullet}}(\Omega^*,\Omega^*)\to
CH_{k}^{Y_{\bullet}}(\Omega^*,\Omega^*)$ which are once again denoted
$d_i$ and $s_j$ by abuse of notation. Using these, the differential 
$D:CH_{\bullet}^{Y_{\bullet}}(\Omega^*,\Omega^*)\to 
CH_{\bullet}^{Y_{\bullet}}(\Omega^*,\Omega^*)$ is defined by setting 
$D(a_0\otimes a_1\otimes\ldots\otimes a_{y_k})$ equal to
$$\sum_{i=0}^{y_k}(-1)^{k+\epsilon_i} 
a_0\otimes\ldots\otimes d_ia_i\otimes\ldots\otimes a_{y_k}+
\sum_{i=0}^k(-1)^id_i(a_0\otimes\ldots\otimes a_{y_k}),$$
where $\epsilon_i$ is again a Koszul sign (see the explicit formula in 
\cite{GTZ}). The simplicial relations imply 
that $D^2=0$ (this is the instance where one uses that $\Omega^*$ is graded 
{\it commutative}). 
These definitions extend by inductive limit to arbitrary (i.e. not 
necessarily finite) simplicial sets.

The homology of $CH_{\bullet}^{Y_{\bullet}}(\Omega^*,\Omega^*)$ w.r.t. the 
differential $D$ is by definition the {\it higher Hochschild homology}
$HH_{\bullet}^{Y_{\bullet}}(\Omega^*,\Omega^*)$ of $\Omega^*$ associated to 
the simplicial set $Y_{\bullet}$. In fact, for two simplicial sets 
$Y_{\bullet}$ and $Y_{\bullet}'$ which have homeomorphic geometric realization,
the complexes $(CH_{\bullet}^{Y_{\bullet}}(\Omega^*,\Omega^*),D)$ and 
$(CH_{\bullet}^{Y_{\bullet}'}(\Omega^*,\Omega^*),D)$ are quasiisomorphic, thus 
the higher Hochschild homology does only depend on the topological space
which is the realization of $Y_{\bullet}$. Therefore we will for example
write $HH_{\bullet}^{\T}(\Omega^*,\Omega^*)$ for the higher Hochschild 
homology of $\Omega^*$ associated to the $2$-dimensional torus $\T$, 
inferring that it is computed w.r.t. some simplicial set having $\T$ as its
geometric realization.

For the simplicial model of the circle $S^1$ given in Example 2.3.1 in 
\cite{GTZ}, one obtains the usual Hochschild homology. In this sense, 
$HH_{\bullet}^{Y_{\bullet}}$ generalizes ordinary Hochschild homology. 
  
Example 2.4.5 in \cite{GTZ} gives:

\begin{prop}
For the simplicial model of the $2$-torus $\T$ given in Example 2.3.2 of
{\it loc. cit.}, the algebra $CH_{\bullet}^{\T}(\Omega^*,\Omega^*)$ is 
quasiisomorphic to 
$$CH_{\bullet}(CH_*(\Omega^*,\Omega^*),CH_*(\Omega^*,\Omega^*)).$$
\end{prop}

In this sense, the holonomy cycle $P({\mathcal A}_0)$ (constructed in the 
previous section) may be regarded as living in 
the higher Hochschild complex $CH_{\bullet}^{\T}(\Omega^*,\Omega^*)$.
This completes the proof of Theorem \ref{main_theorem}.

\begin{rem}
Observe that the element $P({\mathcal A}_0)$ in 
$CH_{\bullet}^{\T}(\Omega^*,\Omega^*)$ is of total degree zero. Recall from 
\cite{GTZ} (Corollary 2.4.7) the iterated 
integral map $It^{Y_{\bullet}}$ of \cite{GTZ} which provides a morphism
of differential graded algebras
$$It^{Y_{\bullet}}:CH_{\bullet}^{Y_{\bullet}}(\Omega^*,\Omega^*)\to
\Omega^{\bullet}(M^{\T},U{\mathfrak h}).$$
The image of $P({\mathcal A}_0)$ in $\Omega^{\bullet}(M^{\T},U{\mathfrak h})$ 
represents a degree
zero cohomology class which associates to each map $f:\T\to M$ an element of 
$U{\mathfrak h}$
which is interpreted as the gerbe holonomy taken over $f(\T)\subset M$.
We believe that an explicit expression of this cohomology class (in the 
special case of an abelian gerbe where all forms are real-valued) is given 
exactly by Gawedski-Reis' formula (2.14) \cite{GawRei}.
The factors $g_{ijk}$ do not appear in our formula, because we did not do the 
gluing yet and therefore everything is local. 

Observe further that following the steps in the proof of Corollary 2.4.4 of 
\cite{GTZ}, one may express $P({\mathcal A}_0)$ in terms of matrices in $A$ 
and $B$. 
\end{rem}

\end{document}